\documentclass[12pt]{article}

\usepackage{a4wide}

\usepackage{amsmath}
\usepackage{amssymb}
\usepackage{amsthm}

\newtheorem{theorem}{Theorem}
\newtheorem{lemma}[theorem]{Lemma}
\newtheorem{proposition}[theorem]{Proposition}

\newtheorem{definition}[theorem]{Definition}
\newtheorem{remark}[theorem]{Remark}

\newcommand{\glie}{\mathfrak{g}}
\newcommand{\hlie}{\mathfrak{h}}
\newcommand{\klie}{\mathfrak{k}}
\newcommand{\llie}{\mathfrak{l}}
\newcommand{\nlie}{\mathfrak{n}}
\newcommand{\plie}{\mathfrak{p}}
\newcommand{\qlie}{\mathfrak{q}}
\newcommand{\tlie}{\mathfrak{t}}
\newcommand{\slie}{\mathfrak{s}}
\newcommand{\ulie}{\mathfrak{u}}

\newcommand{\N}{\mathbb{N}}

\newcommand{\R}{\mathbb{R}}
\newcommand{\C}{\mathbb{C}}

\newcommand{\dolbeault}{\overline{\partial}}
\newcommand{\lnil}{[\![}
\newcommand{\rnil}{]\!]}

\title{Maximal hypoellipticity and Dolbeault cohomology
	representations for $\mathrm{U}(p,q)$}
\author{Nicolas Prudhon, IRMA Strasbourg}
\date{\today}

\begin{document}

\maketitle
\abstract{
Let $Y=G/L$ be a flag manifold for a reductive $G$ and $K$ a maximal 
compact subgroup of $G$.
We define here an equivariant differential operator on $G/L\cap K$ 
playing the role of an equivariant Dolbeault Laplacian for the complex manifold
$G/L$, using a distribution transverse  to the fibers of 
$G/L\cap K \to G/L$ and satisfying the H\"ormander condition. 
We prove here that this operator is not maximal hypoelliptic.
}

\section*{Introdution}
There  are two challenging problems in representation theory of Lie groups. 
The first one is to classify unitary representations for large classes of 
Lie groups. 
Connected nilpotent Lie groups form such a class, and Kirillov
established, for any connected nilpotent Lie group,
a bijective correspondance between the set of coadjoint
orbits and the set of (equivalence classes of) unitary irreducible representations
of the group. This approach lead to the second problem : to realize 
unitary representations geometrically. These two problems are still open 
for reductive groups, but the technique of coadjoint orbits is a constant source
of inspiration.  For reductive groups there are three kind of orbits \cite{vogan}: 
the hyperbolic orbits,
the elliptic orbits, and the nilpotent orbits. The hyperbolic orbits lead
to the theory of parabolic induction and Knapp-Stein intertwining operators.
This is appropriate to construct unitary representations that are weakly contained
in the regular representation.
The  elliptic orbits are related with the theory of cohomological
induction and the geometry of flag manifolds. The study of nilpotent
orbits lead to the theory of unipotent representations. We are concerned here with
the geometry of flag manifolds and we use the theory of coadjoint orbits for nilpotent
Lie groups to handle regularity problems of differential operators on flag manifolds.

Let $G$ be a reductive Lie group and $Y$ be a flag manifold for $G$.
The $G$-space $Y$ is a complex manifold with an equivariant
complex structure, and is a homogeneous space of the form $G/L$,
where the Lie subgroup $L$ is reductive but don't need to be compact.
A representation $\chi$ of $L$ is chosen,
and the usual Dolbeault complex is twisted by $\chi$. The smooth 
cohomology $H^*(\dolbeault_\chi)$ of this complex is proved by  H. W. Wong \cite{wong:95}
to be a Fr\'echet representation of $G$ whose underlying Harish-Chandra module
is isomorphic to the cohomologically induced representation $R(\chi)$.
The proof of H. W. Wong uses the double fibration
$G/L \leftarrow G/L\cap K \rightarrow G/K$, where the group $K$
is a maximal compact subgroup of $G$.  
One conjecture
that if $\chi$ is a unipotent unitary representation of $L$, whatever it means, then the 
representation $H^*(\dolbeault_\chi)$  is unitarizable. However, as a Fr\'echet
space it can not carry a unitary structure. In the best case, when $L$ is compact,
one choose a $G$-invariant hermitian metric on $Y$ and then consider
two objects : the Hilbert space of square integrable sections of the twisted
Dolbeault complex, and the Dolbeault laplacian 
$\overline{\square}_\chi=\dolbeault_\chi \dolbeault_\chi^* + \dolbeault_\chi^* \dolbeault_\chi$.
This differential operator is elliptic and is a selfadjoint operator on the Hilbert
space. Its $L^2$-kernel is then proved to be a unitary representation that 
infinitesimally isomorphic to the Fr\'echet representation. Such representations are
sums of discete series \cite{atiyah-schmid},\cite{connes-moscovici:82}. 
In the general case, necessary to find other representations,
the flag manifold does not carry any $G$-invariant hermitian metric. 
A positive metric is then defined in \cite{rsw} to define the Hilbert space,
and I proved in full generality \cite{prudhon:06} that this Hilbert space is a continuous 
$G$-module. The proof again uses the double fibration considered by Wong.
To replace the $G$-invariant selfadjoint operator, an invariant non-positive form on $Y$ is
defined \cite{rsw},\cite{bkz}. It is used to define the adjoint 
$\dolbeault_{\chi,\mathrm{inv}}^*$ and 
the harmonic space $\ker \dolbeault_\chi \cap \ker \dolbeault_\chi^*$. The point
is that the invariant operator 
$\dolbeault_\chi\dolbeault_{\chi,\mathrm{inv}}^*+
\dolbeault_{\chi,\mathrm{inv}}^*\dolbeault_\chi$  
does not satisfy any regularity condition such as ellipticity and can not be used. 

We propose
here a new invariant operator, defined via the fibration $\pi_L \colon G/L\cap K \to G/L$ and study its
regularity properties as an operator on $G/L\cap K$.  We first define a distribution $E$
transverse to the fibers that satisfies the H\"ormander's condition. It is used to pullback
the Dolbeault operator also denoted by $\dolbeault$. 
The manifold $G/L\cap K$ has a $G$-invariant positive
metric defined by the Killing form, and we can use it to define the formal
adjoint $\dolbeault^*$ of the pullback of the Dolbeault operator. We then define
$\overline{\square}= \dolbeault\dolbeault^*
+\dolbeault^*\dolbeault $. We first show that on sections constant along
the fibers, this operator equals (up to an operator of lower order)
the H\"ormander Laplacian which is known to be maximal hypoelliptic. We next show
that on the whole space of sections over $G/L\cap K$ the operator 
$\overline{\square}$ is not maximal hypoelliptic.  To prove this
we provide  the tangent space of  $G/L\cap K$ with a nilpotent algebra structure, 
canonically associated to the fibration $\pi_L$,
and find a non trivial irreducible representation $\sigma$ of the associated connected nilpotent
Lie group such that the image by $\sigma$ of the $E$-symbol of 
$\overline{\square}$ is not injective on the space of smooth
vectors of $\sigma$. Actually this turns out to be the case for many representations.

The representation
$\chi$ would have been of interest for the (more delicate) questions of positivity 
for instance but
does not come into questions of regularity ; we then use the usual 
Dolbeault complex.

During the preparation of this article I benefited of many helpfull discussions  with L. Barchini, P. Julg, J.M. Lescure, B. Nourrigat, R. Ponge and A. Valette.

\section{The Dolbeault Laplacian}
\subsection{Definition}
Let $Y=G/L$ be a flag manifold for a reductive 
Lie group $G$. This means that $Y$ is an open orbit in a flag manifold
$G^\C / Q$, where $G^\C$ is the complexified Lie group of $G$
and $Q$ is a parabolic subgroup of $G^\C$ ; we also require that
$Y$ admits a $G$-invariant measure. We note $\glie_0$ the Lie algebra
of $G$, and $\glie$ its complexification and use the same convention
with other real and complex Lie algebras and spaces.
Then $Y$ has an equivariant complex structure.
Choices of a maximal compact subgroup $K$ of $G$ and of a base point $y_0 \in Y$
can be made such that the reductive group $L=\mathrm{Stab}_G(y_0)$
is the centralizer of a compact torus with Lie algebra $\tlie'_0 \in \glie_0$,
$L^\C$ is the Levi part of $Q$ and $K/L \cap K$ is a maximal compact complex
submanifold of $Y$. The parabolic algebra $\qlie$ has a decompostion
$\qlie = \llie \oplus \ulie$, and $\glie = \llie \oplus \ulie \oplus \overline{\ulie}$	
with $X \mapsto \overline X$ is the conjugaison associated to the real form 
$G$ of $G^\C$. The space $\ulie$ is $L$-isomorphic to the antiholomorphic
tangent space $T^{0,1}_e G/L$.
Note that the connected reductive subgroup $L$ need not to be compact,
so that $Y$ does not have a $G$-invariant Riemannian metric in general. 

The manifold $Y$ has a $G$-invariant complex structure : this means that
the De Rham operator $d$ writes $d=\partial + \dolbeault$, where
$\partial \colon \wedge^{p,q}TY_\C \to \wedge^{p+1,q}TY_\C$ and
$\dolbeault \colon \wedge^{p,q}TY_\C \to \wedge^{p,q+1}TY_\C$ are
$G$-equivariant operators.
The restriction to $\wedge^{0,*}TY=\wedge^*\ulie$ 
of the operator $\dolbeault$ is called the Dolbeault operator . 
The manifold $Z=G/L\cap K$ fibers over $Y$ and the group $G$ acts on 
it properly. It then admits a $G$-invariant Riemannian metric.
We define the horizontal space at a point $z\in Z$ to be the orthocomplement $E_z$
of the space $F_z$ tangent at $z$ to the fiber trough $z$.

We then have a connexion $E$ on the fibration $\pi_L$ which enables to pullback 
the Dolbeault operator.
\begin{proposition}
	Let $Y$ be a complex manifold with $G$-invariant complex structure
	and $\pi \colon Z \to Y$ an equivariant fibration, with fiber $F$.
	We suppose that the exact sequence
	$$
	TF \to TZ \to \pi^*TY
	$$
	has an equivariant splitting. Let $p_*^{0,1}$ be the transposed map $p_*$
	of this splitting followed by the projection to the (pullback of the)
	antiholomorphic tangent space $\pi^*T^{0,1}Y$. Then there exists a unique 
	operator $\dolbeault'$ on $Z$ satisfying the following conditions.
	\begin{align}
		\dolbeault'(\pi^*\omega) &=\pi^*(\dolbeault\omega) \\
		\dolbeault'(f\pi^*\omega) &=p_*^{0,1}df\wedge(\pi^*\omega)
				+f\pi^*(\dolbeault\omega)
	\end{align}
\end{proposition}
The operator $\dolbeault'$ will be denoted $\dolbeault$ when no confusion arises.
\begin{proof}
	We have to check that, for any $f\in C^\infty(Z)$, any $g\in C^\infty(Y)$ non zero,
	and any $\omega\in \Gamma(Y, \wedge^{*}T^{0,1}Y)$, we have 
	$\dolbeault'(f\pi^*g\pi^*(g^{-1}\omega))=\dolbeault'(f\pi^*\omega)$. Now,
	\begin{align*}
		\dolbeault'(f\pi^*g\pi^*(g^{-1}\omega)) 
		&= p_*^{0,1}d(f\pi^*g)\wedge\pi^*(g^{-1}\omega) 
			+ f\pi^*g\pi^*(\dolbeault(g^{-1}\omega)) \\
		&= p_*^{0,1}(df)\pi^*g\pi^*(g^{-1})\wedge\pi^*(\omega)
			+ fp_*^{0,1}d(\pi^*g)\pi^*(g^{-1})\wedge\pi^*(\omega) \\
		&\phantom{=} {} 
			+ f\pi^*g\pi^*\dolbeault(g^{-1})\pi^*\omega
			+ f\pi^*(g)\pi^*(g^{-1})\pi^*(\dolbeault\omega) \\
		&= p_*^{0,1}df\wedge(\pi^*\omega) +f\pi^*(\dolbeault\omega)
			+ f\pi^*(g^{-1}\dolbeault g 
			+ g\dolbeault(g^{-1}))\wedge\omega \\
		&= \dolbeault'(f\pi^*\omega)\,.
	\end{align*}
\end{proof}
The action of $G$ on $G/L\cap K$ is proper. In particular, this $G$-space admits
a $G$-invariant Riemannian metric. As usual the choice of such a metric enables to 
define a bilinear pairing $(\,,\,)$ between the space of forms with compact supports 
and the space of forms.
The $*$-operator is then given by $(\alpha,\beta)d\mathrm{vol}=\alpha\wedge(*\beta)$.
We then define the adjoint of the pullbacked Dolbeault operator (on homogeneous forms) by
$$ 
\dolbeault^*\omega = (-1)^{|\omega|}(*\dolbeault *)\omega\,.
$$
It remains to define the Dolbeault Laplacian by 
$$
\overline{\square}=\dolbeault\dolbeault^*+\dolbeault^*\dolbeault\,.
$$
This operator is $G$-equivariant by construction.
The question is : can we build an algebra of pseudodifferential operators
on which the Dolbeault Laplacian admits a parametrix~? The result we prove here
gives a negative answer to that question.

\subsection{Structure of the transerve subbundle}
The Dolbeault Laplacian is clearly not elliptic. To study more involved regularity
properties of this operator, we will need detailed information on the bundle $E$.
We now investigate the structure of this bundle. 

\begin{definition} 
	A subbundle $E$ of the tangent space $TZ$ is a 2-step {\em bracket generating}
	subbundle if for any point $p \in Z$, the space $[X,Y](p) \mod E_{p}$, with $X$ and $Y$
	running over sections of $E$,is the whole space $T_pZ / E_p$. 
	In particular the bundle homomorphism
	\begin{equation} \label{hormander2}
	[\,,\,]_0 \colon \bigwedge\!{}^2 E \longrightarrow TZ/E\,,
	\end{equation}
	induced by the barcket $[\,,\,]$ of vectors fields, is onto.
	We say that $E$ satisfies
	the {\em H\"ormander condition at order $2$}.
\end{definition} 
\begin{lemma} \label{existence}
	The subbundle $E$ of the tangent bundle satisfies the H\"ormander
	condition at order $2$.
\end{lemma}
\begin{proof}
	Without lost of generality we may assume that $\glie$ is simple.
	It is enough to prove that 
	$\slie=\ulie \oplus \overline{\ulie} + [\ulie,\overline{\ulie}]$
	is a non zero ideal of $\glie$. As $\qlie$ is a parabolic subalgebra, $\ulie$
	and $\overline{\ulie}$ are sums of root spaces. Moreover $\glie=\llie \oplus \ulie 
	\oplus \overline{\ulie}$, so we get
	$$
	[\ulie,\overline{\ulie}]=
	([\ulie,\overline{\ulie}]\cap \llie) \oplus
	([\ulie,\overline{\ulie}]\cap \ulie) \oplus
	([\ulie,\overline{\ulie}]\cap \overline{\ulie})\,.
	$$
	Using this one checks that 
	$[\ulie\oplus\overline\ulie,[\ulie,\overline{\ulie}]] \subset \slie$.
	Let $X,X' \in [\ulie,\overline{\ulie}]$.
	We write $X=X_l + X_u + X_{\overline{u}}$ thanks to
	the preceding equation, and $X'=[X'_u,X'_{\overline{u}}]$.
	This gives : $[X,X']=[[X_l,X'_u],X'_{\overline{u}}]+
	[X'_u,[X_l,X'_{\overline{u}}]]+X^{''}$ with $X^{''}\in \slie$.
	So $[X,X'] \in \slie$ and $\slie$ is a subalgebra of $\glie$. It is also
	clearly stable by $\llie$.
\end{proof}
We now state a more precise result when $G$ is the group $U(p,q)$ and 
$L=U(p_1)\times U(p_2,q)$ with $p_1+p_2=p$.
\begin{lemma} \label{unique}
	There exists a sequence $\Gamma=(\gamma_1,\ldots,\gamma_r)$ of roots in 
	$\Delta(\llie\cap \plie)$, such that, for any $\alpha\in\Delta(\ulie)$
	there exists at most one $1\leq i \leq r$ and $\beta \in \Delta(\ulie)$
	such that $\alpha \pm \gamma_i=\beta$. Moreover, such an $\alpha$ exists
	for all compact roots in $\Delta(\ulie)$ or it exists for all non compact roots
	in $\Delta(\ulie)$.
\end{lemma}
\begin{proof}
	The roots in $\Delta$ are $e_i-e_j$ and
	\begin{align*}
		\Delta(\llie\cap\plie)=\{e_i-e_j\,; p_1 < i\leq p < j \leq p+q\}\\ 
		\Delta(\ulie)=\{e_i-e_j\,; 1\leq i\leq p_1 < j \leq p+q\}	\,.
	\end{align*}
	Set $r=\mathop{min}\{p_1,q\}$ (the real rank of the noncompact semi simple
	part of $\llie$) and let $\Gamma=(\gamma_i)$ be any set of strongly orthogonal 
	roots.
	For exemple, one may take $\gamma_i=e_{p_1+i}-e_{p+q-i}$. The result follows 		
	easily. In fact, if $\alpha=e_i-e_j \in \Delta(\ulie)$ then the only
	$\beta=e_k-e_l$ that may work are those with $k=i$ or $l=j$, and only one
	of them can lies in $\Gamma$.
\end{proof}
\begin{remark}
	This lemma is also easily seen to be true when $G$ is any real rank $1$ group.
\end{remark}

\subsection{Statement of the main result} 
Let us precise now the regularity property of differential operators we
want to investigate. Let $X_1, \ldots ,X_k$ be vector fields on a neighborhood $V$ of a point $x_0 \in \R^n$, and
let $E_{x_0}$ be the subspace of $\R^n$ generated by the vectors $X_i(x_0)$.
We also assume that vectors $[X_i,X_j](x_0) \mod (E_{x_0})$ generate the vector space $\R^n / E_{x_0}$.
The space of operators of order less than $m$ is
the space of operator $P$ that can be written in the form
\begin{equation}\label{E-operators}
	P= \sum_{|\alpha|\leq m}a_\alpha(x) X^\alpha\,,\quad 
	X^\alpha=X_1^{\alpha_1}\cdots X_k^{\alpha_k},
\end{equation}
where the coefficient $a_\alpha(x)$ are smooth functions of the variable $x$ on $V$.
\begin{definition} \cite{helffer-nourrigat}
	A differential operator $P$ of order $m$ is maximal hypoelliptic
	at $x_0$, if there exists a neighborhood $V$ of $x_0$ and a constant 
	$C$ so that for all $u\in C^\infty_c(V)$,
	$$
	\sum_{|\alpha|\leq m} \|X^\alpha u\|_{L^2}
	\leq C\big( \|u\|_{L^2}+\|Pu\|_{L^2} \big) \,.
	$$
\end{definition}
Maximal hypoellipticity of an operator $P$ implies that $P$ is hypoelliptic, i.e. 
$$
Pu \text{ smooth} \quad \Rightarrow  \quad u \text{ smooth}\,.
$$
The principal $E$-symbol is by definition $p=\sum_{|\alpha| = m}a_\alpha(x) \xi^\alpha$.
We will use the sign "$\simeq$" to say that two operators have the same
principal $E$-symbol.	
The following result is well known (see \cite{helffer-nourrigat}).
\begin{proposition}
	The H\"ormander Laplacian $\sum_i X_i^2$ is maximal hypoelliptic.
\end{proposition}

We choose here the metric given by the Killing form $B$.
More precisely, the metric is defined at the origin by 
$$
\langle X,Y \rangle= -B(X,\theta(\overline{Y}))\,.
$$
This form is definite positive on $\glie$ and is $\klie$-invariant. The tangent 
spaces $T_eY\simeq \ulie\oplus\overline{\ulie}$ and 
$T_eZ \simeq \ulie\oplus\overline{\ulie}\oplus (\llie\cap\plie)$ are provided
with this hermitian metric.
	
\begin{theorem} \label{main:theorem}
	Let $G=\mathrm{U}(p,q)$ and $L=\mathrm{U}(p_1)\times\mathrm{U}(p_2,q)$, 
	with $p_1+p_2=p$. 
	The Laplacian $\overline{\square}$ is not maximal hypoelliptic 
	at the origin $eL\cap K$.
\end{theorem}
One may conjecture this this result is true for any semisimple Lie group
and flag manifold.
The exposition of the proof is intended to make clear that only
the lemma \ref{unique} as to be generalized. So let $G$ be a reductive
Lie group with a compact Cartan subalgebra, and $G/L$ be a flag manifold for $G$.
This assumption on the Cartan subalgebra makes less technical
the computation of the principal $E$-symbol , but we should
proceed without it.

The next section is devoted to the proof of the theorem \ref{main:theorem}.
To prepare the proof we compute here the local expression of the principal 
$E$-symbol of this operator. The Cartan subalgebra being compact, 
we may suppose that 
$$
\tlie_0 \subset 
\llie_0 \cap \klie_0 \subset 
\llie_0 \subset 
\glie_0 \,.
$$
Let $\Delta$ be the root system of the pair $(\glie,\tlie)$
All roots of the root system $\Delta(\glie,\tlie)$ 
being compact or non compact, it makes sens to define $\Delta(\ulie\cap\klie)$
and $\Delta(\ulie\cap\plie)$ and so on.
We choose a system $\Delta^+(\glie,\tlie)$
of positive roots such that $\Delta(\ulie)\subset\Delta^+(\glie,\tlie)$.
As the Killing
form is non-degenerate there exists for any $\alpha \in \Delta$
a vector $H_\alpha \in \tlie$ so that for all $H\in \tlie$, 
$\alpha(H)=B(H,H_\alpha)$.
\begin{lemma}
	There exists  an orthonormal basis $(E_\alpha)_{\alpha\in\Delta}$ of
	root vectors satisfying
	\begin{subequations}\label{eq:helgason}
		\begin{align}
			[E_\alpha,E_{-\alpha}] &= H_\alpha 
			\label{eq:helgason:a}\\
			[E_\alpha,E_\beta] = N_{\alpha,\beta}E_{\alpha+\beta} 
			\quad &\text{with }N_{\alpha,\beta}=0 
			\text{ if }\alpha + \beta \notin \Delta  
			\label{eq:helgason:b}\\
			N_{\alpha,\beta} &= -N_{-\alpha,-\beta}\,.
			 \label{eq:helgason:c}
		\end{align}
	\end{subequations}
\end{lemma}
\begin{proof}
	According to \cite[theorem 5.5]{book:helgason:62} 
	there exists a basis $(E'_\alpha)$
	satisfying equations 
	(\ref{eq:helgason}). The relation (\ref{eq:helgason:a}) implies that 
	$B(E'_\alpha,E'_{-\alpha})=1$. Moreover $B(E'_\alpha,E'_\beta)=0$ if 
	$\alpha+\beta\neq 0$, and $\|E'_\alpha\|>0$, so it follows that
	$-\theta(\overline{E'_\alpha})=c_{-\alpha}E'_{-\alpha}$, 
	with $c_\alpha c_{-\alpha}=1$.
	We now define $E_\alpha=x_\alpha E'_\alpha$ where $x_\alpha x_{-\alpha}=1$ 
	and $x_\alpha^2=-c_\alpha$. 
	We then get
	\begin{gather}
		-\theta(\overline{E_\alpha})=x_\alpha c_{-\alpha}E'_{-\alpha}
			=-E_{-\alpha}\quad\text{and} \label{involution} \\
		[E_\alpha,E_{-\alpha}]=x_\alpha x_{-\alpha}[E'_\alpha,E'_{-\alpha}]
			=H_\alpha \,. \label{brackets}
	\end{gather}
	So that $\langle E_\alpha,E_\alpha\rangle = B(E_\alpha,E_{-\alpha}) = 1$ and the
	basis $(E_\alpha)$ is now orthonormal. Using equations (\ref{involution}) 
	and (\ref{brackets}) is easy to check that the basis $(E_\alpha)$
	again satifies the equations (\ref{eq:helgason}).
\end{proof}
We set $\overline{Z}_\alpha=E_\alpha$ and
$$
Z_\alpha=\left\{
\begin{array}{ll}
	-E_{-\alpha} & \text{if $\alpha$ is compact,} \\
	E_\alpha & \text{if $\alpha$ is non compact.}
\end{array} 
\right.
$$
This notation is concording with the complex structure.
Let us now define the real vectors $X_\gamma$ and $Y_\gamma$. 
$$
X_\gamma=\frac{1}{\sqrt{2}}\big(Z_\gamma+\overline{Z}_\gamma \big)\,, \quad
Y_\gamma=-\frac{i}{\sqrt{2}}\big(\overline{Z}_\gamma -Z_\gamma \big)\,.
$$
The sytem 
$(X_\gamma,Y_\gamma)_{\gamma \in \Delta^+\setminus \Delta^+(\mathfrak{l}\cap \mathfrak{k})}$
 is an orthonormal basis  of $T_eZ$ and 
$(X_\gamma,Y_\gamma)_{\gamma \in \Delta(\mathfrak{u})}$ is an orthonormal 
basis of $E_e \simeq T_eY$.
Moreover if $J$ denotes the complex multiplication operator, one has
$Y_\gamma=JX_\gamma$, for $\gamma\in \Delta(\ulie)$ (and $\gamma \in \Delta(\llie \cap \plie)$
when $L/L\cap K$ is a hermitian symertic space).
We also have
\begin{subequations}\label{real-complex}\begin{align}
	X_\alpha &= \frac{1}{\sqrt{2}}(E_\alpha-E_{-\alpha})  
	&Y_\alpha &=-\frac{i}{\sqrt{2}}(E_\alpha+E_{-\alpha})
	\quad    &\text{if $\alpha$ is compact,} \label{real-complex-1}\\
	X_\beta  &= \frac{1}{\sqrt{2}}(E_\beta+E_{-\beta})
	&Y_\beta  &=-\frac{i}{\sqrt{2}}(E_\beta-E_{-\beta})
	\quad    &\text{if $\beta$ is non compact.} \label{real-complex-2}
\end{align}\end{subequations}
\begin{proposition} \label{field:structure}
	For $\alpha \in \Delta(\mathfrak{u}\cap \mathfrak{k})$ and 
	$\beta\in\Delta(\mathfrak{u}\cap\mathfrak{p})$ we have 
	\begin{subequations}\label{lie_structure}\begin{align}
		[X_\alpha,X_\beta] &= \frac{1}{\sqrt{2}}
		\Big(
		N_{\alpha\beta}
		X_{\alpha+\beta}+N_{\alpha,-\beta}X_{|\alpha-\beta|}
		\Big)\\
		[X_\alpha,Y_\beta] &= \frac{1}{\sqrt{2}}
		\Big(
		N_{\alpha,\beta}Y_{\alpha+\beta}
		-\epsilon(\alpha-\beta)N_{\alpha,-\beta}Y_{|\alpha-\beta|}
		\Big)\\
		[Y_\alpha,X_\beta] &= \frac{1}{\sqrt{2}}
		\Big(
		N_{\alpha,\beta}Y_{\alpha+\beta}
		+\epsilon(\alpha-\beta)N_{\alpha,-\beta}Y_{|\alpha-\beta|}
		\Big)\\
		[Y_\alpha,Y_\beta] &= -\frac{1}{\sqrt{2}}
		\Big(
		N_{\alpha,\beta}X_{\alpha+\beta}
		-N_{\alpha,-\beta}X_{|\alpha-\beta|}
		\Big)
	\end{align}\end{subequations}
	The vectors involving roots of the form $\alpha+\beta$ lie in $E_e$.
	The vectors involving roots of the form $\alpha-\beta$ may lie in $F_e$,
	but don't need to.  Other brackets of base vectors lie in $E_e$. 
\end{proposition}
To prove this proposition one just computes using 
equations (\ref{eq:helgason:b},\ref{eq:helgason:c})
and the fact that 
if $\alpha\in\Delta(\ulie\cap\klie)$ is compact and $\beta\in\Delta(\ulie\cap\plie)$ 
is non compact then $\alpha\pm\beta$ either is a non compact root or is not a root.

Let $e_\gamma$ be the exterior multiplication by $Z_\gamma$. Then 
the Dolbeault operator has the following principal $E$-symbol.
$$
\dolbeault\simeq\sum_{\gamma\in\Delta(\ulie)} e_\gamma \overline{Z}_\gamma \,,
$$
where $\overline{Z}_\gamma$ is here the left invariant vector field generated
by $\overline{Z}_\gamma$. Let $i_\gamma$ be the interior multiplication
by $\overline{Z_\gamma}$ with respect to the chosen metric. Then
$$
\dolbeault^*\simeq-\sum_{\gamma\in\Delta(\ulie)} i_\gamma Z_\gamma\,.
$$
According to the previous notations these equations become
\begin{gather*}
	\dolbeault\simeq \sum_{\gamma\in\Delta(\ulie)} \frac{e_\gamma}{\sqrt{2}}
		\big(X_\gamma-iY_\gamma \big) \,, \\
	\dolbeault^*\simeq-\sum_{\gamma\in\Delta(\ulie)} \frac{e_\gamma}{\sqrt{2}}
		\big(X_\gamma+iY_\gamma \big)\,.
\end{gather*}
It now remains to compute.
$$
\overline{\square}
\simeq -\frac{1}{2}\Big( 
\sum_{\gamma} e_\gamma\left(X_\gamma-iY_\gamma\right)\cdot
\sum_{\gamma'} i_{\gamma'}\left(X_{\gamma'}+iY_{\gamma'}\right) +
\sum_{\gamma'} i_{\gamma'}\left(X_{\gamma'}+iY_{\gamma'}\right)\cdot
\sum_{\gamma} e_\gamma\left(X_\gamma-iY_\gamma\right)
\Big) \,.
$$
Let us write the diagonal terms separetely.
\begin{align*}
	\overline{\square} \simeq
	&{-\frac{1}{2}} 
		\sum_{\gamma\in\Delta(\ulie)}
		\left(e_\gamma i_\gamma+i_\gamma e_\gamma \right)
		\left(X_\gamma^2+Y_\gamma^2 \right) \\
	&{-\frac{1}{2}} 
		\sum_{\gamma\neq\gamma'} e_\gamma i_{\gamma'}
		\Big[
		\left( X_\gamma X_{\gamma'}+Y_\gamma Y_{\gamma'}\right)
		+i\left(X_\gamma Y_{\gamma'}+Y_\gamma X_{\gamma'}\right)
		\Big]\\
	&{}\phantom{-\frac{1}{2} 
		\sum_{\gamma\neq\gamma'}} +  i_{\gamma'}e_\gamma
		\Big[
		\left( X_{\gamma'}X_\gamma+Y_{\gamma'}Y_\gamma\right)
		+i\left( Y_{\gamma'}X_\gamma+ X_{\gamma'}Y_\gamma\right)
		\Big]
\end{align*}
We have $e_\gamma i_{\gamma'}+i_{\gamma'}e_\gamma = \delta_{\gamma\gamma'}$
(Kronecker' symbol). 
\begin{align*}
	\overline{\square} \simeq
	&{-\frac{1}{2}} \sum_{\gamma\in\Delta(\ulie)}
		\left(X_\gamma^2+Y_\gamma^2 \right) \\
	&{-\frac{1}{2}} \sum_{\gamma\neq\gamma'} e_\gamma i_{\gamma'}
		\Big[
		\left( [X_\gamma,X_{\gamma'}]+[Y_\gamma,Y_{\gamma'}] \right)
		+i\left( [X_\gamma,Y_{\gamma'}]+[Y_\gamma,X_{\gamma'}]\right)
		\Big]
\end{align*}
Using proposition \ref{field:structure} one gets
\begin{gather*}
	[X_\gamma,X_{\gamma'}]+[Y_\gamma,Y_{\gamma'}]=
		\sqrt{2}\,N_{\alpha,-\beta}X_{|\alpha-\beta|}  \quad\text{and}\\
	[X_\gamma,Y_{\gamma'}]+[Y_\gamma,X_{\gamma'}]=
		-\sqrt{2}\,N_{\alpha,-\beta}Y_{|\alpha-\beta|}\,,
\end{gather*}
if $\gamma=\alpha$ is compact and $\gamma'=\beta$ is non compact. One
has similar relations when $\gamma=\beta$ is non compact and $\gamma'=\alpha$ is
compact. Other brackets are horizontal and they don't appear in
the principal $E$-symbol. This gives
\begin{align}\label{local}
	\overline{\square} \simeq
	&-\frac{1}{2} 
		\sum_{\gamma\in\Delta(\ulie)}
		\left( X_\gamma^2+Y_\gamma^2 \right)
		+\frac{\sqrt{2}}{2} \sum_{\gamma\in\Delta(\llie\cap\plie)}
		\Big[ \left(\sum \!{}^* N_{\alpha,\-\beta} 
		\left(e_\alpha i_\beta - e_\beta i_\alpha\right)\right)X_\gamma \notag \\
	&{}\phantom{-\frac{1}{2} 
		\sum_{\gamma\in\Delta(\ulie)}
		\left( X_\gamma^2+Y_\gamma^2 \right)
			+\frac{\sqrt{2}}{2} 
		\sum_{\gamma\in\Delta(\llie\cap\plie)}
		\Big[ }
		+i \left(\sum \!{}^* N_{\alpha,\-\beta} 
		\left(e_\alpha i_\beta + e_\beta i_\alpha\right)\right)Y_\gamma 
		\Big]
\end{align}
where the sums $\sum^*$ are over $\alpha\in\Delta(\ulie\cap\klie)$, 
$\beta\in\Delta(\ulie\cap\plie)$ and $|\alpha-\beta|=\gamma$. The 
local formula (\ref{local}) will be used later in the proof of theorem 
\ref{main:theorem}.

This formula is already usefull for functions. In fact, the terms of
classical order $1$ vanish on functions, so $\overline\square$ is
maximally hypoelliptic when restricted to functions  because
it has the same principal $E$-symbol as the H\"ormander Laplacian (up to
a constant).

%
\section{The Rockland condition}

\subsection{Hypoellipticity criterion} 
For the proof of the theorem \ref{main:theorem} 
we use techniques of Folland and  Stein \cite{folland-stein}.
We now provide the tangent space $T_eZ$ with a nilpotent Lie algebra structure $\nlie_0$.
This structure is given by the brackets $[\,,\,]_0$, and the identification
of $TZ/E$ with $F$.
The Lie brackets $\lnil \,,\, \rnil$ is then given as follows. Compare with
proposition \ref{field:structure}.
\begin{definition}\label{lie:structure}
For $\alpha \in \Delta(\mathfrak{u}\cap \mathfrak{k})$ and 
	$\beta\in\Delta(\mathfrak{u}\cap\mathfrak{p})$ we have 
	\begin{subequations}\begin{align}
		\lnil X_\alpha,X_\beta\rnil &= \frac{1}{\sqrt{2}}
		\Big(
		N'_{\alpha,-\beta}X_{|\alpha-\beta|}
		\Big)\\
		\lnil X_\alpha,Y_\beta\rnil &= \frac{1}{\sqrt{2}}
		\Big(
		-\epsilon(\alpha-\beta)N'_{\alpha,-\beta}Y_{|\alpha-\beta|}
		\Big)\\
		\lnil Y_\alpha,X_\beta\rnil &= \frac{1}{\sqrt{2}}
		\epsilon(\alpha-\beta)N'_{\alpha,-\beta}Y_{|\alpha-\beta|}
		\Big)\\
		\lnil Y_\alpha,Y_\beta\rnil &= -\frac{1}{\sqrt{2}}
		\Big(
		-N'_{\alpha,-\beta}X_{|\alpha-\beta|}
		\Big),
	\end{align}\end{subequations}
	where $N'_{\alpha,-\beta}=N_{\alpha,-\beta}$ if 
	$\alpha-\beta\in \Delta(\llie\cap\plie)$ and $0$ otherwise. 
	All other brackets of base vectors are defined to be $0$.
\end{definition}

Let $P$ be a differential operator on an open set of $\R^n$ as in the first part, 
with principal $E$-symbol $p$. We say that $P$ satisfies the Rockland condition
if for any unitary irreducible non trivial representation
$\pi$ of the simply connected nilpotent Lie group 
$N=\exp(\nlie_0)$, the operator $\pi(p)$ is injective
on the space of smooth vectors of $\pi$.
The sympbol $p$ is seen here as an element of the enveloping algebra $U(\nlie)$ of
$\nlie$.

\begin{theorem} \cite{helffer-nourrigat}
	The following are equivalent
	\begin{enumerate}
		\item $P$ has a parametrix in the $E$-pseudodifferential calculus,
		\item $P$ satisfies to the Rockland condition,
		\item $P$ is maximal hypoelliptic.
	\end{enumerate}
\end{theorem}

Let $N$ be a nilpotent Lie group with Lie algebra $\nlie_0$. Then $N$
acts on $\nlie^*_0$ by the coadjoint representation. Kirillov defined a one-to-one
correspondance between coadjoint orbits and (equivalence classes of) 
irreducible unitary representations of the group $N$ constructed in three 
steps as follows.
\begin{lemma} 
	Let $l$ be a form on $\nlie_0$ and 
	$B_l \colon (X,Y) \mapsto l([X,Y])$. Then there exists
	an isotropic subalgebra $\hlie_0$ of $\nlie_0$ for $B_l$
	such that $\mathrm{codim}\hlie_0=\frac{1}{2}\mathrm{rank}B_l$.
\end{lemma}
\noindent Then $\exp(il)$ is a one dimensionnal representation of the nilpotent 
group $H=\exp(\hlie_0)$.
\begin{lemma}
	The induced representation $\mathrm{Ind}_H^N e^{il}$ is irreducible
	and its class only depends on the coadjoint orbit of $l$.
\end{lemma}
\noindent There is also a converse statement.
\begin{lemma} 
	All the irreducible representations of $N$ arise in this way exactly once.
\end{lemma}
 
We will need to recognize induced representations realized on $\R^n$.
Let $\pi$ be a representation of the nilpotent Lie algebra $\nlie_0$ on 
$\mathcal{S}(\R^n)$.
We suppose that, for any $X\in \nlie_0$ , the operator $\pi(X)$ has the form
$$
\pi(X) = \sum_{k=1}^{n-1}P_k(y_1,\ldots,y_{k-1};X)\frac{\partial}{\partial y_k}
		    +iQ(y_1,\ldots y_n ; X)\,,
$$
where $P_k(\cdot;X)$ and $Q(\cdot ;X)$ 
are polynomials on $\R^n$ depending linearily on $X$. We also assume
that the linear forms $\xi_k(X)=P_k(0;X)$ are linearily independent.
Let $l$ be the linear form on $\nlie_0$ defined by
$l(X)=Q(0;X)$ and $\hlie_0= \cap \ker \xi_k$.
\begin{proposition}\cite[Proposition 1.6.1]{helffer-nourrigat}\label{rep-induite}
	Under the above assumptions, the subspace $\hlie_0$ is a subalgebra 
	of $\nlie_0$, isotropic for $B_l$. Moreover, the representation $\pi$
	is unitarily equivalent to $\mathrm{Ind}_H^Ne^{il}$.
\end{proposition}

\subsection{Proof of the main theorem}
Here we  prove that the evaluation of the Dolbeault laplacian has a kernel of positive
dimension on many representations under conditions on root systems. The choice
of these representations and the proofs of the root systems conditions are made
for the groups $G=\mathrm{U}(p,q)$  and $L=\mathrm{U}(p_1)\times\mathrm{U}(p_2,q)$,
with $p_1+p_2=p$. One can expect that this can be done in full generality.

We now have to find unitary irreducible representations
of the connected nilpotent Lie group $H$ and to realize them on $L^2(\R^n)$.
This will lead to a partial differential equation on $\R^n$. In other words the
linear form on $\nlie_0$ that gives the representation of $N$,  has to be taken 
such that the obtained partial differential
equation (can be solved and) has a non zero solution space.
 Let $l \in \nlie_0^*$ be a linear form on $\nlie_0$ with coordinates $(\xi_\gamma,\eta_\gamma)$
in the dual basis of $(X_\gamma,Y_\gamma)$. Let $\pi_l$ be the representation
of $N$ associated to the coadjoint orbit of $l$. 

Using definiton \ref{lie:structure} one find that 
the form $B_l \colon (X,Y) \mapsto l([X,Y])$ as a martix of the form
$$
\left( \begin{array}{ccc} 0& A & 0 \\ -A^\text{t} & 0 & 0 \\ 0&0&0 \end{array} \right)
\,, \quad \text{with } 
A=\frac{N'_{\alpha,-\beta}}{\sqrt{2}}\left( \begin{array}{cc} 
	\xi_{|\alpha-\beta|}& 
	-\varepsilon(\alpha-\beta)\eta_{|\alpha-\beta|} \\ 
	\varepsilon(\alpha-\beta)\eta_{|\alpha-\beta|} &
	\xi_{|\alpha-\beta|}
\end{array} \right)_{\alpha,\beta}\,.
$$
We make the following assumption on $l$.

\noindent
(H) \hspace{1,2cm}$A$ has a maximal rank.

\noindent
If hypothesis (H) is true then

\noindent
(H') \hspace{1cm}
either $\plie_0$ or $\llie_0\cap \plie_0 \oplus (\ulie\oplus\overline\ulie)\cap\klie_0$ is
a maximal abelian subalgebra of $\nlie_0$. 

\noindent 
This means that the hypothesis (H) is more an hypothesis on the pair $(G,Q)$ than
on the linear form $l$. Let us assume hypothesis (H). Let 
$\hlie_0$ be the abelian subalgebra of $\nlie_0$ such that 
$$
\hlie_0=\plie_0\quad\text{if}\quad
\dim \plie_0 =\mathrm{max}\big\{\dim \plie_0\,;\, \dim \llie_0\cap \plie_0 \oplus (\ulie\oplus\overline\ulie)\cap\klie_0\big\}\,,
$$
and $\hlie_0=\llie_0\cap \plie_0 \oplus (\ulie\oplus\overline\ulie)\cap\klie_0$ otherwise. Then
$$
\mathrm{codim} \hlie_0 = \frac{1}{2} \mathrm{rank} B_l\,,
$$
and $\hlie_0$ is an isotropic subspace for $B_l$.
\begin{lemma}\label{H}
	Let $G=\mathrm{U}(p,q)$ and $L=\mathrm{U}(p_1)\times\mathrm{U}(p_2,q)$.
	There exists a linear form $l$ such that hypothesis (H) is satisfied.
\end{lemma}
\begin{proof}
	Take $l$ be non zero on root vectors correszponding
	 to a set of strongly orthogonal roots in $\Delta^+(\llie\cap\plie)$ such
	as in the proof of lemma \ref{unique}, and
	$0$ elsewhere. Then $A$ is "diagonal" with no zero on the diagonal,
	by lemmas \ref{existence},\ref{unique}. 
\end{proof}

\noindent\emph{First case}. Let us begin with the case $\hlie_0=\plie_0$.
Let $s=\mathrm{dim}_\C K/L\cap K=\mathrm{dim}\ulie\cap\klie$.
Then $\pi_l=\mathop{Ind}_H^N e^{il}$ 
is a unitary irreducible representation of $N$ on $L^2(\nlie_0/\hlie_0)$ that 
can be seen as a representation on $L^2(\R^{2s})$.
We note $(x_\alpha,y_\alpha)_{\alpha\in\Delta(\ulie\cap\klie)}$ the
canonical basis of $\R^{2s}$. Thanks to proposition \ref{rep-induite}, we have
\begin{gather*}
	\pi_l(X_\alpha)=\frac{\partial}{\partial x_\alpha} + i\xi_\alpha\,,
	\quad
	\pi_l(Y_\alpha)=\frac{\partial}{\partial y_\alpha} + i\eta_\alpha\,, \\
	\pi_l(X_\beta)=i\sum_{\alpha}
		\left[
		\frac{N'_{\alpha,-\beta}}{\sqrt{2}} 
		\left(
		\xi_{|\alpha-\beta|}x_\alpha -\varepsilon(\alpha-\beta)\eta_{|\alpha-\beta|}y_\alpha
		\right)
		\right]
		+ i\xi_\beta \,,\\
	\pi_l(Y_\beta)=i\sum_{\alpha}
		\left[
		\frac{N'_{\alpha,-\beta}}{\sqrt{2}} 
		\left(
		\varepsilon(\alpha-\beta)\eta_{|\alpha-\beta|}x_\alpha+\xi_{|\alpha-\beta|}y_\alpha
		\right) 
		\right]
	+ i\eta_\beta\,,\\
	\pi_l(X_\gamma)=i\xi_\gamma\,, 
	\quad
	\pi_l(Y_\gamma)=i\eta_\gamma \,.
\end{gather*}
To make the computation more easy we also suppose
that 
\begin{equation} \label{orbite}
	\xi_\alpha=\eta_\alpha=\xi_\beta=\eta_\beta=0\,.	
\end{equation}
This is a priori not true in general that any
orbits admits a form of this kind, but this is enough, to prove the theorem, to
find such forms such that $\pi_l(\overline\square)$ is not injective.
Then, the operator $\pi_l\left(\overline\square\right)$ has the following form. 
\begin{equation}\label{rep-lap}
\pi_l\left(\overline\square\right)=-\frac{1}{2}\sum_\alpha
\left[\frac{\partial^2}{\partial x_\alpha^2} +  \frac{\partial^2}{\partial y_\alpha^2}
-r_\alpha^2(x_\alpha^2 +y_\alpha^2)\right]
+\sum_\alpha \left[
\sum\!^{*}
M_{\alpha,\beta} \right] \,,
\end{equation}
where $r_\alpha$ is the positive real number such that 
$r_\alpha^2 = \sum\!^* \frac{N'_{\alpha,-\beta}\,^2}{2} 
(\xi_{|\alpha-\beta|}^2+\eta_{|\alpha-\beta|}^2)$, and
$$
M_{\alpha,\beta}=\frac{iN'_{\alpha,-\beta}}{\sqrt{2}}
\Big[
(\xi_{|\alpha-\beta|}+i\eta_{|\alpha-\beta|})e_\alpha i_\beta 
-(\xi_{|\alpha-\beta|}-i\eta_{|\alpha-\beta|})i_\beta e_\alpha
\Big] 
$$
is an endomorphism of $\wedge^* \ulie$ and the sum $\sum\!^*$
is over the set of roots $\beta\in\Delta(\ulie\cap\plie)$ such that
$\alpha-\beta \in \Delta(\llie\cap\plie)$.

Let $D_\alpha=-\frac{1}{2}\left[\frac{\partial^2}{\partial x_\alpha^2} +  
\frac{\partial^2}{\partial y_\alpha^2}
-r_\alpha^2(x_\alpha^2 +y_\alpha^2)\right]$ and $M_\alpha=\sum\!^*M_{\alpha,\beta}$.
We have to find eigenvalues of $\sum_\alpha D_\alpha$ and 
$\sum_\alpha M_\alpha$ of opposite signs
and the same absolute value. Making the change of variables 
$$
x_\alpha \mapsto r_\alpha^{\frac{1}{2}}x_\alpha \,\quad
y_\alpha \mapsto r_\alpha^{\frac{1}{2}}y_\alpha \,,
$$
the operator $D_\alpha$ becomes 
$-\frac{r_\alpha}{2}\left[\frac{\partial^2}{\partial x_\alpha^2} +  
\frac{\partial^2}{\partial y_\alpha^2}
-(x_\alpha^2 +y_\alpha^2)\right]$. It is $-\frac{r_\alpha}{2}$ 
times the Hermite operator of dimension $2$.
Its eigenvalues are then $-kr_\alpha$, with $k\in\N^*$. As the operators
$D_\alpha$ differentiate on different variables, we see that  the eigenvalues
of $\sum_\alpha D_\alpha$ are $-\sum_\alpha k_\alpha r_\alpha$, with
$k_\alpha \in \N^*$. We also note that the eigenfunctions of the Hermite
operator are of the form $P(x)e^{-\frac{\|x\|^2}{2}}$ where $P$ is a polynomial
in $x=(x_1,\ldots,x_{2s})$. So they are in the
Schwarz space, so are smooth vectors of the representation $\pi_l$.

Let us now show that $\pm\sum_\alpha r_\alpha$ is an eigenvalue 
of $\sum_\alpha M_\alpha$.  We first show that $r_\alpha$ is 
an eigenvalue of $M_\alpha$. Let $\Delta(\ulie\cap\klie)=\{\alpha_1,\ldots,\alpha_s\}$ 
and $v=Z_{\alpha_1}\wedge\cdots \wedge Z_{\alpha_s}$. If $\beta\neq\beta'$, then
$M_{\alpha,\beta}M_{\alpha,\beta'}(v)=M_{\alpha,\beta'}M_{\alpha,\beta}(v)=0$
and moreover 
$$
M_{\alpha,\beta}^2(v)= \frac{N^{'2}_{\alpha,-\beta}}{2}
(\xi_{|\alpha-\beta|}^2+\eta_{|\alpha-\beta|}^2)\,.
$$
It follows that
$$
M_\alpha^2(v)=\sum\!^* M_{\alpha,\beta}^2(v)=r_\alpha^2 v\,.
$$
So the vector $\pm r_\alpha v + M_\alpha v$ is an eigenvector for $M_\alpha$
with eigenvalue $\pm r_\alpha$.
\begin{proposition}
	Let $k\leq s$ and $\{i_1;\cdots ;i_k\}\subset \{1;\cdots ;s\}$. Then 
	$\prod_{l=1}^k M_{\alpha_{i_l}}v$ does not depend
	on the order of the $i_l$.
\end{proposition}
This proposition is easily checked by induction on $k$.
We now define by induction, for $k\leq s$, the vectors $v_k$ by $v_0=v$ and
$$
v_k = (r_{\alpha_k}+M_{\alpha_k})v_{k-1}\,.
$$ 
The preceding proposition shows that if $v_{k-1}$
is an eigenvector for $M_{\alpha_l}$, $l < k$, with 
eigenvalue $r_{\alpha_l}$, then $v_k$ is again an eigenvector
for $M_{\alpha_l}$, $l < k$, with 
eigenvalue $r_{\alpha_l}$. 
\begin{lemma}
	Let $G=\mathrm{U}(p,q)$ and $L=\mathrm{U}(p_1)\times\mathrm{U}(p_2,q)$.
	There exists a linear form $l$ on $\nlie_0$ satisfying hypothesis (H), and
	such that $v_k$ is an eigenvector for $M_{\alpha_k}$, with  eigenvalue
	$r_{\alpha_k}$.
\end{lemma}
\begin{proof}
	Defining $l$ has in the proof of lemma \ref{H} again works.
\end{proof}
Hence $v_s$
is a simultaneous eigenvector for all $M_\alpha$'s, 
with respective eigenvalue $r_\alpha$.
So $v_s$ is an eigenvector for $\sum_\alpha M_\alpha$
with eigenvalue $\sum_\alpha r_\alpha$.
We end this first case $\hlie=\plie_0$ remarking that
the constructed eigenvector lies in $\wedge^s \ulie$, and
this means that $\overline\square$ is not maximally
hypoelliptic on degree $s=\dim_\C K/L\cap K$.

\vspace{0.3cm}
\noindent \emph{Second case}. Let us now assume that
$t=\dim \ulie\cap \plie < s$. Switching the role played
in the first case by the $\alpha$'s and the $\beta$'s, one
similiraly proves that $\overline \square$ is not maximal
hypoelliptic on $\wedge^t \ulie$. Using the duality
$$
\wedge \colon \wedge^t \ulie \otimes \wedge^s \ulie 
\rightarrow \wedge^{\mathrm{max}}\ulie\,,
$$
one shows that $\overline \square$ is not 
maximal hypoelliptic on degree $s$ in  the second case too.

Finally, we have shown that $\overline \square$ is never
maximal hypoellitpic on degree $s$ and on the complementary degree $t$.
It would be remarkable if these degrees are the only one where 
this phenomena arises. It is obviously the case for $G=\mathrm{U}(2,2)$
and $L=\mathrm{U}(1)\times \mathrm{U}(1,2)$ for instance, because any coadjoint orbit
admits a linear form satisfying equation (\ref{orbite}).

\bibliographystyle{alpha} 
\bibliography{biblio}


\end{document}